\documentclass[12pt]{amsart}

\usepackage[mathscr]{eucal}
\usepackage{cite}

\newcommand{\dif}{\mathrm{d}}
\DeclareMathAlphabet{\pazocal}{OMS}{zplm}{m}{n}

\usepackage{amssymb}
\input cyracc.def

\newfont{\tricyr}{wncyr10 at 12pt}
\newfont{\tricyi}{wncyi10 at 12pt}
\newfont{\tricyb}{wncyb10 at 12pt}
\newfont{\Tricyr}{wncyr10 at 13.6pt}
\newfont{\Tricyi}{wncyi10 at 13.6pt}
\newfont{\Tricyb}{wncyb10 at 13.6pt}
\newfont{\tricmr}{cmr10 at 13.6pt}
\newfont{\tricmi}{cmti10 at 13.6pt}
\newfont{\tricmb}{cmb10 at 13.6pt}

\usepackage{enumerate}

\usepackage{graphicx}

\makeatletter
\@namedef{subjclassname@2020}{%
  \textup{2010} Mathematics Subject Classification}
\makeatother

\newtheorem{theorem}{Theorem}[section]

\newtheorem{lemma}[theorem]{Lemma}

\theoremstyle{definition}
\newtheorem{definition}[theorem]{Definition}
\newtheorem{Not}{Remark}

\numberwithin{equation}{section}

\begin{document}

\baselineskip=17pt



\title[Existence and non-existence of solutions]{Nontrivial solutions to the Dirichlet problems for semilinear degenerate elliptic equations}

\author[D. T. Luyen]{Duong Trong Luyen}
\address{International Center for Research and Postgraduate Training in Mathematics, Institute of Mathematics,Vietnam Academy of Science and Technology,Hanoi,Vietnam}
\address{Department of Mathematics\\ Hoa Lu University\\ Ninh Nhat\\ 
Ninh Binh city, Vietnam.}
\email{dtluyen.dnb@moet.edu.vn, dtluyen@math.ac.vn}

\author[N. M. Tri]{Nguyen Minh Tri}
\address{Institute of Mathematics\\ Vietnam Academy of Science and Technology \\ 18 Hoang Quoc Viet, 10307
 Cau Giay, Hanoi, Vietnam.}
\email{triminh@math.ac.vn}

\author[D. A. Tuan]{Dang Anh Tuan}
\address{University of Sciences\\ Vietnam National University\\ 
334 Nguyen Trai, Thanh Xuan, Hanoi, Vietnam.}
\email{datuan1105@gmail.com}
\date{}

\begin{abstract}
 In this article, we study the existence of non-trivial weak solutions for the following boundary-value problem 
 \begin{gather*}
 -\frac{\partial^2 u}{\partial x^2} -\left|x\right|^{2k}\frac{\partial^2 u}{\partial y^2}=f(x,y,u) \quad\text{ in }\Omega, \
 u=0  \quad\text{ on }\partial\Omega,
 \end{gather*}
where $\Omega$ is a bounded domain with smooth boundary in $\mathbb{R}^2, \Omega \cap \{x=0\}\ne \emptyset,$ $k >0,$ $f(x,y,0)=0. $

\end{abstract}



\maketitle

\section{Introduction}
The existence and nonexistence of nontrivial solutions to the Dirichlet boundary value problems of semilinear degenerate elliptic differential equations have been extensively studied for more than 30 years, see for example \cite{D. Jerison:1981a}, \cite{D. Jerison:1987a}, \cite{N. M. Tri:1998a}, \cite{N. M. Tri:1998b}, \cite{thuy:2012}, \cite{Luyen:2015}, \cite{Hua:2020b}, \cite{N. M. Tri:2010a}, \cite{N. M. Tri:2014a}, \cite{Nga:2023a}  and the references therein. The methods used here are modeled from the classical elliptic theory. To deal with the nonexistence of nontrival solutions we need to establish a type of the Pohozaev identity for the problems under consideration. To obtain the exitence of nontrivial solutions we can use the variational method that is first we construct corresponding weghted Sobolev spaces, then prove embedding theorems for them, and finaly use the critical point theory. The Yamabe problem on the Heizenberg group  and CR manifolds were considered in \cite{D. Jerison:1981a}, \cite{D. Jerison:1987a}. The boundary value problems for semilinear degenerate elliptic equations  equations involving the Grushin operator was treated in \cite{N. M. Tri:1998a}, \cite{N. M. Tri:1998b}. Then the results of \cite{N. M. Tri:1998a}, \cite{N. M. Tri:1998b} were generalized for semilinear  equations containing the strongly degenerate elliptic operator in \cite{thuy:2012}. The non-linear terms of equations in \cite{N. M. Tri:1998a}, \cite{N. M. Tri:1998b}, \cite{thuy:2012} do not depend on the space variable.  Semilinear equations involving the $\Delta_\gamma$-Laplace operator was considered in \cite{Lanconelli:2012}. The class of  $\Delta_\gamma$-Laplace operators containing both the Grushin and strongly degenerate elliptic operators. The multiplicity of solutions for semilinear $\Delta_\gamma$ equations was considered in \cite{Luyen:2015}. Estimates of Dirichlet eigenvalues for degenerate $\Delta_\gamma$-Laplace operator was given in \cite{Hua:2020b}. Semilinear equations involving the Grushin operator with critical exponent was treated in \cite{Huong:2023a}. Nontrivial solutions for degenerate elliptic equations in a solid torus was investigated in a recent paper \cite{Nga:2023a}.

In this paper, we are mainly concerned with the semilinear subelliptic Dirichlet problem
\begin{alignat}{2}
 -\frac{\partial^2 u}{\partial x^2} -\left|x\right|^{2k}\frac{\partial^2 u}{\partial y^2}&=f(x,y,u)&\qquad& \text{in}\quad \Omega,\label{Luyen2} \\
u&=0 &\qquad &\text{on}\quad \partial \Omega, \label{Luyen3}
\end{alignat}
where $\Omega$ is a bounded domain with smooth boundary in $\mathbb{R}^2, \Omega \cap \{x=0\}\ne \emptyset, k >0. $ We emphasize that the nonlinear term of this equation  depends on the space variable. We will see that the critical exponet for this problem essentially depends on the space variable. The critical exponent here is higher than the one obtained in \cite{N. M. Tri:1998a}, \cite{N. M. Tri:1998b}. This is a familiar phenomenon happened in the Henon equation with respect to the semilinear Laplace equation, see \cite{Ni:1982a}. 

The structure of our note is as follows:  In Section 2, we establish a nonexistence theorem via an identity of
Pohozaev’s type. In Section 3, we prove some embedding theorems for weighted
Sobolev spaces associated with the problem. In Section 4, we establish the existence of weak solutions.

\section{Nonexistence results}

In this section, we deal with the nonexistence of nontrivial solutions to the problem \eqref{Luyen2}-\eqref{Luyen3} for $f(x,y,\xi)=\left|x\right|^{2k}\left|\xi\right|^{p-1}\xi, p\ge 1$. Put
\begin{align*}
F(x,y,\xi):=\int\limits_0^\xi f(x,y,\tau)\dif \tau =\frac{\left|x\right|^{2k}}{p+1}|\xi|^{p+1}, 
\end{align*}
and let $\nu =(\nu_{x},\nu_{y})$ be the unit outward normal on $\partial\Omega$. 
We will denote by $\mathscr S^2(\overline\Omega)$ the linear space of functions $u\in C(\Omega)$
such that
\begin{align*}
\frac{\partial u}{\partial x}, \left|x\right|^k\frac{\partial u}{\partial y}, 
\frac{\partial^2 u}{\partial x^2}, \left|x\right|^{2k}\frac{\partial^2 u}{\partial y^2}
\end{align*}
exist in the weak sense of distributions in $\Omega$ and can be extended to $\overline\Omega$.\\
A function $u(x,y)\in \mathscr S^2(\overline\Omega)$ is said to be a solution to the problem
 \eqref{Luyen2}-\eqref{Luyen3} if 
\begin{alignat}{2}
 -\frac{\partial^2 u}{\partial x^2} -\left|x\right|^{2k}\frac{\partial^2 u}{\partial y^2}&=\left|x\right|^{2k}\left|u\right|^{p-1}u&\qquad& \text{in}\quad \Omega,\notag \\
u&=0 &\qquad &\text{on}\quad \partial \Omega. \notag
\end{alignat}
\begin{lemma}\label{bd1}
 Let $u(x,y)\in \mathscr S^2(\overline\Omega)$ be a
solution of the problem \eqref{Luyen2}-\eqref{Luyen3}. Then we have
\begin{align}\label{1.3}
\left(\frac{2+3k}{p+1}-\frac{k}{2}\right)\int\limits_{\Omega}\left|x\right|^{2k}\left|u\right|^{p-1}u^2 \dif x\dif y
= \frac{1}{2} \int\limits_{\partial\Omega} [x\nu_x+(1+k)y\nu_y](\nu^2_x+\left|x\right|^{2k}\nu^2_y) \left(\frac{\partial u}{\partial \nu}\right)^2
\dif s.
\end{align}
\end{lemma}
\begin{proof}
The proof of this lemma is similar to that of  Lemma 1 in \cite{thuy:2012}. We omit the details.
\end{proof}
\begin{definition}
A domain $\Omega$ is called $G_k-$
starshape with respect to the origin  if $(0,0)\in \Omega$ and 
$x\nu_x+(1+k)y\nu_y \ge 0$ at every point of $\partial\Omega.$
\end{definition}
From Lemma \ref{bd1}, we have that
\begin{theorem}\label{DL1}
 Let $\Omega$ be $G_k-$starshaped with respect to the origin 
and 
$$ 
p>\frac{4+5k}{k}.
 $$
Then the problem \eqref{Luyen2}-\eqref{Luyen3} has no nontrivial solution $u\in \mathscr S^2(\overline\Omega)$. 
\end{theorem}
\begin{Not}
From Theorem  \ref{DL1} we can suppose that the number $\frac{4+5k}{k}$
is the critical exponent for the problem \eqref{Luyen2}-\eqref{Luyen3}.
Without the condition of $G_k-$starness of the domain $\Omega$
Theorem \ref{DL1}  may not be true. 
\end{Not}

\section{Embedding theorem}
In this section, we prove an embedding theorem of the Sobolev type.
\begin{definition}
Let $\Omega$  be a bounded domain in $ \mathbb R^2$. 
By $L_k^p(\Omega) \ (1 \le p < +\infty)$ we will denote the set of all measurable functions
$u:\Omega\to \mathbb R$ such that
$$ 
\int\limits_\Omega \left|x\right|^{2k}\left|u\right|^p\dif x\dif y <\infty.
 $$
We define the norm in $L_k^p(\Omega)$ as follows
$$
\left\| u \right\|_{L_k^p(\Omega)}:
=\left(\int\limits_\Omega \left|x\right|^{2k}\left|u\right|^p\dif x\dif y\right)^{\frac{1}{p}}.
$$
By $\mathscr S^{1}_2(\Omega) $ we will denote the set of all functions
 $ u \in L^2(\Omega) $ such that 
$ \frac{\partial u}{\partial x}, \left|x\right|^k \frac{\partial u}{\partial y} \in L^2(\Omega).$
We define the norm in $\mathscr S^{1 }_2(\Omega)$ as follows
$$
\left\|u\right\|_{\mathscr S^{1}_2(\Omega)}:=
\left\{\int\limits_\Omega u^2\dif x\dif y 
+\int\limits_\Omega \left|\nabla_{G_k}u\right|^2\dif x\dif y \right\}^{\frac12},
$$
where $\nabla_{G_k}u:=\left(\frac{\partial u}{\partial x}, \left|x\right|^k \frac{\partial u}{\partial y}\right)$.

We can also define the scalar product in $ {\mathscr S^{1}_2(\Omega)} $ as follows
$$
\left(u , v\right)_{\mathscr S_2^1(\Omega)} 
= \left(u , v\right)_{L^2(\Omega)} + 
\left(\nabla_{G_k} u,\nabla_{G_k} v \right)_{L^2(\Omega)}.
$$
The space $\mathscr S_{2,0 }^1(\Omega)$ is defined as the closure of  $C_0^1(\Omega) $ in the space $\mathscr S^{1 }_2(\Omega)$.
\end{definition}
\begin{theorem}\label{DL:3}
Let $\Omega$  be a bounded domain in $  \mathbb R^2$ with smooth boundary $\partial\Omega$ such that 
$(0,0)\in \Omega $. Then the embedding 
\begin{align*}
\mathscr S_{2,0}^1(\Omega) \hookrightarrow L^q_k( \Omega),  \mbox{ where } 1 \le q \le \frac{4+6k}{k}
\end{align*}
is continuous, i.e. there exists a constant $C_q>0$ such that
\begin{align*}
\left\|u\right\|_{L_k^q(\Omega)}\le C_q \left\|u\right\|_{\mathscr S_{2,0 }^1(\Omega)}, \quad \forall u \in \mathscr S_{2,0 }^1(\Omega).
\end{align*}
Moreover,  the embedding 
\begin{align*}
\mathscr S_{2,0 }^1(\Omega)\hookrightarrow\hookrightarrow L^q_k( \Omega),  \mbox{ where } 1 \le q < 2_k:=\frac{4+6k}{k}
\end{align*}
is compact.
\end{theorem}
\begin{proof}
For the continuous embedding it comes from our Sobolev type inequality (\cite{NTT:2022}):
\begin{equation}\label{NTT1}
||u||_{L^q_k(\Omega)}\le C||\nabla_{G_k}u||_{L^2(\Omega)}, 1\le q\le 2_k.
\end{equation}
Using the H\"{o}lder inequality we have
\begin{equation}\label{NTT2}
||u||_{L^q_k(\Omega)}\le ||u||_{L^1_k(\Omega)}^\mu ||u||_{L^{2_k}_k(\Omega)}^{1-\mu}, 1<q<2_k,
\end{equation}
where $\mu=(2_k-q)/(2_kq-q).$ From \eqref{NTT1}-\eqref{NTT2}, in order to prove the compact embedding for $1\le q<2_k$ we need only prove the compact embbeding 
\begin{equation}\label{compact}
\mathscr S_{2,0 }^1(\Omega)\hookrightarrow\hookrightarrow L^1_k( \Omega).
\end{equation}

\noindent We can consider $u\in \mathscr S_{2,0 }^1(\Omega)$ as $u\in \mathscr S_{2,0 }^1(\mathbb{R}^2)$ by setting $u=0$ outside $\Omega$. Similarly for $L^1_k(\Omega)$. Let $\mathcal{F}$ be a bounded subset in $\mathscr S_{2,0 }^1(\mathbb{R}^2)$. In order to prove the compact embedding \eqref{compact}, we will show that  $\mathcal{F}$ is relatively compact subset in $L^1_k(\Omega)$ or  $\mathcal{G}=\left\lbrace s^\ell u:u\in \mathcal{F}\right\rbrace $ is relatively compact in $L^1(\Omega)$ (or $L^1(\mathbb{R}^2)$).
Because  $\mathcal{F}$ is bounded in $L^1_\ell(\mathbb{R}^2)$, $\mathcal{G}$ is bounded in $L^1(\mathbb{R}^2)$. Thus, in order to prove $\mathcal{G}$ is relatively compact in $L^1(\mathbb{R}^2)$, according to Frechet-Kolmogorov, we need only prove that
$$\mathop {\rm sup}\limits_{v\in \mathcal{G}}\int_{\mathbb{R}^2}\left| v\left( x+h_1, y+h_2\right) -v\left(x, y\right)\right| dxdy\to 0 \text{ as } h=(h_1, h_2)\to 0.$$
 Since $\Omega$ is bounded, there is an $a>0$ such that $\Omega\subset (-a, a)\times(-a, a)$. Let $\epsilon \in (0, a/2),$ $\left| h\right|_\infty=\max\{|h_1|, |h_2|\} <\epsilon$, we have
\begin{align*}
\int_{\mathbb{R}^2}\left| v\left( x+h_1, y+h_2 \right) -v\left(x, y\right)\right| dxdy= &\int_{\left| (x, y)\right|_\infty >a}\left| v\left( x+h_1, y+h_2\right) -v\left(x, y\right)\right| dxdy+\\
+&\int_{2\epsilon<\left|x \right| <a, |y|<a}\left| v\left( x+h_1, y+h_2\right) -v\left(x, y\right)\right| dxdy\\
+&\int_{\left| x\right| <2\epsilon, |y|<a}\left| v\left( x+h_1, y+h_2\right) -v\left(x, y\right)\right| dxdy\\
:=& I_1+I_2+I_3.
\end{align*}
Note that supp$ v\subset \Omega\subset (-a, a)\times(-a, a)$ and $|h_j|<\epsilon$, it is not difficult to obtain
\begin{equation}\label{Tequ8}
I_1\le C||u||_{L^{2_k}_k(\Omega)}\epsilon^{1-1/2_k}, I_3 \le C||u||_{L^{2_k}_k(\Omega)}\epsilon^{(2k+1)(1-1/2_k)}.
\end{equation}
Note that $ v\left( x+h_1, y+h_2\right) - v\left(x, y\right)=\int_0^1h\cdot\nabla v\left(x+th_1, y+th_2\right)dt$ 
and $\nabla v=|x|^{2k}\nabla u+\left( 2k x^{2k-1}u\right) (1, 0)$ so
 
\begin{eqnarray*}
I_2&=&\int_{2\epsilon<\left| x\right| <a, |y|<a}\left| \int_0^1	h\cdot\nabla v\left(x+th_1, y+th_2\right)dt\right| dxdy\\
&\le&\left| h\right| \int_0^1\left( \iint_{2\epsilon<\left| x\right| <a, |y|<a}\left|  \nabla v\left(x+th_1, y+th_2\right)\right| dxdy\right) dt\\
&\le&\epsilon \int_{\epsilon<\left| x\right| <a, |y|<a}\left|  \nabla v\left(x, y\right)\right| dxdy
\end{eqnarray*}
\begin{equation}\label{Tequ10}
\le C\epsilon \int_{\epsilon<\left| x\right| <a, |y|<a} (2k|x|^{2k-1} \left|u\right|+ |x|^{2k}|u_x|+x^{2k} \left| u_y\right|)dxdy:= C\epsilon\left[ J_1+J_2+J_3\right]. 	
\end{equation}
Note that $2_k=(6k+4)/k, k>0$, we have
\begin{equation}\label{Tequ11}
J_1\le C||u||_{L^{2_k}_k(\Omega)}, J_2+J_3\le C||\nabla_{G_k}u||_{L^2(\Omega)}.
\end{equation}
From \eqref{Tequ8}-\eqref{Tequ10}-\eqref{Tequ11} we conclude that $\mathcal{F}$ is relatively compact in $L^1_k(\Omega)$. The proof is complete.
\end{proof}
\begin{Not}
From a result in \cite{N. M. Tri:1998b}, we also have the two norms $\left\|u\right\|_{\mathscr S^{1 }_2(\Omega)}$ and
$$ 
\left| \left\|u\right\|\right|_{\mathscr S_{2,0 }^1(\Omega)}:=\left(\int\limits_\Omega \left|\nabla_{G_k}u\right|^2\dif x\dif y \right)^{\frac12}
 $$
are equivalent in $\mathscr S_{2, 0}^1(\Omega)$.
\end{Not}

\section{Existence results}
From now on we suppose that $f(x,y,\xi)$ has only polynomial growth in $\xi$.
\begin{definition}
A function $u\in  \mathscr S_{2,0 }^1(\Omega) $ 
is called a weak solution of the problem  
 \eqref{Luyen2}-\eqref{Luyen3} if  the identity
$$  
\int\limits_{\Omega}{ \nabla_{G_k} u\cdot  \nabla_{G_k} \varphi \dif x \dif y} 
-\int\limits_{\Omega }{f\left(x,y, u \right)\varphi \dif x \dif y}= 0
$$
is satisfied for every $ \varphi \in   \mathscr S_{2,0 }^1(\Omega).$
\end{definition}
We try to find weak solutions of the problem \eqref{Luyen2}-\eqref{Luyen3} as critical points
of a nonlinear functional. To this end we define the functional $\Phi$ on the space
 $ \mathscr S_{2,0 }^1(\Omega)$ as follows
\begin{align}\label{TriLuyen10:CT14}
\Phi(u) = \frac 12 \int\limits_\Omega \left| \nabla_{G_k} u\right|^2\dif x\dif y -\int\limits_\Omega F(x,y,u)\dif x\dif y.
\end{align}
Using  H\"{o}lder's inequality and Theorem \ref{DL:3}, we can easily obtain
\begin{lemma}\label{TriLuyen5: DLP 1}
Assume that  $f: \Omega \times \mathbb{R} \to \mathbb{R}$ is a Carath\'{e}odory function such that
there exist $p\in (2,2_k)$, 
 $ f_1(x,y)\in L_k^{p_1}(\Omega), f_2(x,y)\in L_k^{p_2}(\Omega),$ 
where 
$  p_1/(p_1-1) <  2_k,  pp_2/(p_2-1)\le  2_k,  p_1 > \max\{1,\frac{2_k p_2}{p_2(p-1) +2_k}\}, p_2> 1 $ 
such that
 $$ \left|f(x,y,\xi)\right| \le \left|x\right|^{2k} (f_1(x,y) +f_2(x,y) \left|\xi\right|^{p-1})  \mbox{ almost everywhere in \ \ } \Omega \times \mathbb{R}.$$
Then $ \Phi_1(u)  \in C^1(\mathscr S_{2,0 }^1(\Omega), \mathbb{R} )$ and 
$$ 
\Phi_1'(u)(v) =  \int\limits_{\Omega}f(x,y,u)v\dif x\dif y
 $$
for all $v \in \mathscr S_{2,0 }^1(\Omega)$, where 
$$  \Phi_1\left(u\right) = \int\limits_\Omega F\left(x,y,u\right) \dif x\dif y,  $$
and  $F(x,y,\xi) = \displaystyle\int\limits_{0}^\xi f(x,y,\tau) \dif \tau$.
\end{lemma}
We assume that  $f: \Omega \times \mathbb{R} \to \mathbb{R}$ is a  Carath\'{e}odory function satisfying
\begin{itemize}
\item[(A1)] $f$  is a real Carath\'{e}odory function on $\Omega \times \mathbb{R}$ such that there exist $q_1\in (2,2_k)$, 
 and  constant $C_0\ge 0$ such that
$$
|f(x,y,\xi)|\leq \left|x\right|^{2k}|(\xi|^{q_1-1} + C_0) \quad \forall (x,y,\xi) \in \Omega\times  \mathbb{R};
$$
\item[(A2)]  there exist  $C\in [0,+\infty)$ and $\psi\in L^{1}_k(\Omega)$ such that
$|f(x,y,\xi)|\leq \left|x\right|^{2k} \psi(x,y)$ for every $(x,y)$ in $\Omega$ and $|\xi|\le C$;

\item[(A3)] there exists  a non-positive function $\varphi$ such that $\int\limits_\Omega \varphi(x,y)\dif x\dif y <\infty$ and
  $ \varphi(x,y)\le \frac{f(x,y,\xi)}{\xi}$ for every $(x,y,\xi)\in \Omega\times  \mathbb{R}$;

\item[(A4)]  $f(x,y,0)=0$ for every $(x,y)$ in $\Omega$ and the following limit holds uniformly for a.e. $(x,y)$ in $\Omega$
\begin{align*}
 \lim\limits_{\xi\to 0}\frac{f(x,y,\xi)}{\left|x\right|^{2k}\xi} =0 \mbox{ and  }\lim\limits_{\xi\to +\infty} \frac{f(x,y,\xi)}{\xi} =+\infty
\end{align*} 

\item[(A5)]  $\frac{f(x,y,\xi)}{\xi}$ is increasing in $\xi \ge C$ and
decreasing in $\xi \le -C$ for every  $(x,y)$ in $\Omega$.
\end{itemize}
Our main result in this section is given in the following theorem.

\begin{theorem}\label{thm1}  Suppose that $f$ satisfies  {\rm (A1)-(A5)}. 
Then the boundary value problem \eqref{Luyen2}-\eqref{Luyen3} has a nontrivial weak
solution.
\end{theorem}

From Theorem \ref{TriLuyen5: DLP 1} and  the fact that $ f$ satisfies  {\rm (A1)}, we have
$\Phi$ is well-defined on $\mathscr S_{2,0 }^1(\Omega)$ and 
$\Phi \in C^1(\mathscr S_{2,0 }^1(\Omega),  \mathbb{R})$ with
$$ 
\langle \Phi'(u), v\rangle =  \int\limits_{\Omega}\nabla_{G_k} u \cdot \nabla_{G_k} v \dif x \dif y
- \int\limits_\Omega f\left(x,y,u(x,y)\right)v(x,y) \dif x\dif y
 $$
for all $v \in  \mathscr S_{2,0 }^1(\Omega).$ 

Hence, the weak solutions of  the problem \eqref{Luyen2}-\eqref{Luyen3} are critical points of the functional $\Phi$. 

\begin{definition}\label{definition} \rm
Let $ \mathbf{B}$ be a real Banach space with its dual space $ \mathbf{B}^*$ and  \linebreak
$\Phi \in C^1( \mathbf{B},  \mathbb{R})$.
For $c \in   \mathbb R$ we say that $\Phi$ satisfies the $(C)_c$ condition if for any sequence $\{u_{n}\}_{n=1}^{+\infty} \subset   \mathbf B$ with
$$ 
\Phi(u_n) \to c\mbox{ and } \left(1+ \left\|u_n\right\|_{ \mathbf{B}}\right)\left\|\Phi'(u_n)\right\|_{ \mathbf{B}^*}\to 0,
 $$
there exists a sub-sequence $\{u_{n_k}\}_{k=1}^{+\infty}$ that converges strongly in $\mathbf {B}$.
\end{definition}

We will use the following version of the Mountain Pass Theorem
\begin{lemma}[see \cite{Cerami:1978,Cerami:1980}]\label{Luyen:BDP}
Let $ \mathbf B$ be a real Banach space and let   $\Phi \in C^1( \mathbf B,  \mathbb R)$  satisfy the $(C)_c$
condition for any $c \in  \mathbb R, \Phi(0)=0$  and
\begin{itemize}
\item[(i)] There exist constants  $\rho, \alpha > 0$ such that  $\Phi(u) \ge \alpha, \forall u \in   \mathbf  B, \left\|u\right\|_{ \mathbf  B} =\rho;$

\item[(ii)]  There exists an  $u_1 \in  \mathbf B, \left\|u_1\right\|_\mathbf B \ge \rho$ such that $\Phi(u_1)\le 0$.
\end{itemize}
Then $ \beta:=\inf\limits_{\lambda \in \Lambda}\max\limits_{0\le t \le 1}\Phi(\lambda(t)) \ge \alpha $ is a critical value of  $\Phi$, where
$$ \Lambda: =\{\lambda \in C\left([0;1],   \mathbf B\right): \lambda(0)=0, \lambda(1)= u_1 \}. $$
\end{lemma}

We prove Theorem \ref{thm1} by verifying that all conditions of Lemma \ref{Luyen:BDP} are satisﬁed. First, we
check the condition (i) in Lemma \ref{Luyen:BDP}.

\begin{lemma} \label{lemma2}
Assume that $f$ satisfies conditions {\rm (A1)} and {\rm (A4)}. Then there exist \linebreak
$\rho, \alpha>0$  such that
$$ 
\Phi(u)\geq\alpha, \ \forall u \in \mathscr S_{2,0 }^1(\Omega), \|u\|_{\mathscr S_{2,0 }^1(\Omega)}=\rho.
 $$
\end{lemma}
\begin{proof}
 Suppose by contradiction that
\begin{align*}
\inf \left\{\Phi(u): u \in \mathscr S_{2,0 }^1(\Omega), \|u\|_{\mathscr S_{2,0 }^1(\Omega)}
 = \frac{1}{n}\right\} \leq 0\quad \forall n\in \mathbb{N}.
\end{align*}
Then, there exists a sequence $\{u_{n}\}_{n=1}^{+\infty}$ in $\mathscr S_{2,0 }^1(\Omega)$ such that
$\|u_{n}\|_{\mathscr S_{2,0 }^1(\Omega)} =\frac{1}{n}$ and \linebreak
 $\Phi(u_{n}) < \frac{1}{n^{3}}$. Hence, we have
\begin{align*}
\frac{1}{n}> \frac{\Phi(u_n)}{\|u_n\|_{\mathscr S_{2,0 }^1(\Omega)} ^2}= 
\frac{1}{2}- \int\limits_{\Omega}\frac{F(x,y,u_{n}(x,y))}{\|u_{n}\|_{\mathscr S_{2,0 }^1(\Omega)} ^{2}}\dif x\dif y,
\end{align*}
and thus
\begin{align}\label{TriLuyen: CT0}
\int\limits_{\Omega}\frac{F(x,y,u_{n}(x,y))}{\|u_{n}\|_{\mathscr S_{2,0 }^1(\Omega)} ^{2}}\dif x\dif y
>  \frac{1}{2}- \frac{1}{n}.
\end{align}
By {\rm (A4)} for each $\varepsilon>0$, we can find a number $\delta>0$ such that
\begin{align} \label{TriLuyen: CT1}
\left| F(x,y,\xi)\right| \le {\varepsilon \left|x\right|^{2k} \xi^2} \mbox{ for  } \left|\xi\right| \le \delta.
\end{align}
From {\rm (A1)}, we deduce that 
\begin{align}\label{TriLuyen: CT2}
\left| F(x,y,\xi)\right| \le \frac{1}{q_1} \left|x\right|^{2k} \left|\xi \right|^{q_1} +C(\delta) \left|x\right|^{2k} \left|\xi \right|^{q_1} \mbox{ for  } \left|\xi\right| \ge \delta.
\end{align}
It follows from \eqref{TriLuyen: CT1}, \eqref{TriLuyen: CT2},  Theorem \ref{DL:3} and 
H\"{o}lder's  inequality that
\begin{align*}
 &\left| \int\limits_{\Omega}F(x,y,u_{n}(x,y))\dif x\dif y\right| \\
\le &\varepsilon\int\limits_{\Omega} { \left|x\right|^{2k} \left|u_n(x,y)\right|^2} \dif x\dif y + 
\frac{1}{q_1}\int\limits_{\Omega} \left|x\right|^{2k} \left|u_n(x,y) \right|^{q_1} \dif x\dif y
+ C(\delta)\int\limits_{\Omega} \left|x\right|^{2k} \left|u_n(x,y) \right|^{q_1}\dif x\dif y \\
\le& \varepsilon C_1 \left\|u_n\right\|_{\mathscr S_{2,0 }^1(\Omega)}^2 +C_2 \left\|u_n\right\|_{\mathscr S_{2,0 }^1(\Omega)}^{q_1} 
+ C(\delta)C_3 \left\|u_n\right\|_{\mathscr S_{2,0 }^1(\Omega)}^{q_1},
\end{align*}
hence 
\begin{align*}
 \left| \int\limits_{\Omega}\frac{F(x,y,u_{n}(x,y))}{\|u_{n}\|_{\mathscr S_{2,0 }^1(\Omega)} ^{2}}\dif x\dif y \right| 
 \to 0 \mbox{ as } n \to +\infty,
\end{align*}
which yields a contradiction to \eqref{TriLuyen: CT0}. Lemma \ref{lemma2} is proved.
\end{proof}
Next, we check the condition (ii) in Lemma \ref{Luyen:BDP}.
 \begin{lemma} \label{lemma1}
Let $\rho$ be as in Lemma \ref{lemma2} and 
assume that $f$ satisfies conditions {\rm (A3)} and {\rm (A4)}. Then
there exists $u_1$ in $\mathscr S_{2,0 }^1(\Omega)\setminus B(0,\rho)$ such that $\Phi(u_1)< 0$.
\end{lemma}
\begin{proof}
Take a point 
$u\in \mathscr S_{2,0 }^1(\Omega)$ such that 
$\|u\|_{\mathscr S_{2,0 }^1(\Omega)}=1, u >0$  and $\|u\|_{L^2(\Omega)}\ne 0.$ Then, for any constant $R>0$
$$ 
\Phi(Ru) =\frac{R^2}{2}- \int\limits_{\Omega}F(x,y, Ru(x,y))\dif x\dif y.
 $$
Since {\rm (A4)}, there exists a number $M>0$ such that $f(x,y,\xi) \ge \frac{4\xi}{\|u\|^2_{L^2(\Omega)}}$ for $\xi\ge M;$ 
moreover, from {\rm (A3)} we get
\begin{gather*}
F(x,y,\xi)\ge \int\limits_{0}^M \varphi(x,y)\tau \dif \tau=\varphi(x,y)\frac{M^2}{2} \mbox{ for } 0\le \xi\le M,\\
F(x,y,\xi) =\int\limits_0^M f(x,y,\tau)\dif\tau+\int\limits_{M}^\xi f(x,y,\tau)\dif\tau \\
\ge \varphi(x,y)\frac{M^2}{2} + \frac{2\xi^2}{\|u\|^2_{L^2(\Omega)}}-\frac{2M^2}{\|u\|^2_{L^2(\Omega)}} \mbox{ for } \xi\ge M.
\end{gather*}
As $0\le Ru\le M$ on $\Omega^u_{\frac{M}{R}}:=\{(x,y)\in \Omega: \left| u(x,y)\right|\le \frac{M}{R}\}$, we have
\begin{gather*}
\int\limits_{\Omega}F(x, y, Ru(x,y))\dif x\dif y = \int\limits_{\Omega^u_{\frac{M}{R}}}F(x,y, Ru(x,y))\dif x\dif y
 + \int\limits_{\Omega\backslash {\Omega^u_{\frac{M}{R}}}}F(x,y, Ru(x,y))\dif x\dif y \\
\ge  \int\limits_{\Omega\backslash {\Omega^u_{\frac{M}{R}}}}\frac{2R^2u^2(x, y)}{\|u\|^2_{L^2(\Omega)}}\dif x\dif y
+ {M^2}\int\limits_{\Omega} \varphi(x,y)\dif x\dif y - \frac{2M^2}{\|u\|^2_{L^2(\Omega)}}\operatorname{meas}(\Omega),
\end{gather*}
where $\operatorname{meas}(\cdot)$ denotes the Lebesgue measure of a set in  $\mathbb{R}^2$.
By a theorem of Lebesgue, there exists a number $R_0$ such that
$$ 
\int\limits_{\Omega\backslash {\Omega^u_{\frac{M}{R_0}}}}u^2(x,y) \dif x\dif y  \ge \frac{\|u\|^2_{L^2(\Omega)}}{2}.
 $$
Therefore, if $R\ge R_0$ then 
$$ 
\int\limits_{\Omega}F(x,y, Ru(x,y))\dif x\dif y \ge R^2 + {M^2}\int\limits_{\Omega} \varphi(x,y)\dif x\dif y - \frac{2M^2}{\|u\|^2_{L^2(\Omega)}}\operatorname{meas}(\Omega).
 $$
Therefore, if 
$$ 
R> \max\left\{2\sqrt{-{M^2}\int\limits_{\Omega} \varphi(x,y)\dif x\dif y + \frac{2M^2}{\|u\|^2_{L^2(\Omega)}}\operatorname{meas}(\Omega)}, R_0\right\},
 $$
then 
$$ 
\Phi(Ru) \le {-{M^2}\int\limits_{\Omega} \varphi(x,y)\dif x\dif y+ \frac{2M^2}{\|u\|^2_{L^2(\Omega)}}\operatorname{meas}(\Omega)} -\frac {R^2}{2}<0.
 $$
Thus, $\Phi$ satisfies the condition (ii) in Lemma \ref{Luyen:BDP}.
\end{proof}

 \begin{lemma} \label{lemmcc0}
Assume that $f$ satisfies conditions {\rm (A2)} and {\rm (A5)}. 
Then  there exists a positive real
number $C_3$ such that
\begin{align*}
 f(x,y,s)s-2F(x,y,s)\le f(x,y,t)t - 2F(x,y,t) +C_3\left|x\right|^{2k}\psi(x,y),\quad \forall (x,y)\in\Omega,|s|\le |t|.
\end{align*}
\end{lemma}
\begin{proof}
Since {\rm (A5)}, by Lemma 2.3 in \cite{Liu:2010}, we have that,  for any $(x,y)\in \Omega,$ the mapping
\begin{align*}
\xi \mapsto f(x,y,\xi)\xi-2F(x,y,\xi)
\end{align*}
is increasing in $\xi \ge C$ and decreasing in $\xi \le -C$. 
Hence 
\begin{align*}
f(x,y,s)s-2F(x,y,s) \le f(x,y,t)t - 2F(x,y,t), \quad \forall (x,y)\in\Omega, C\le s\le t.
\end{align*}
  Let $(x,y)\in\Omega$ and $\xi \in [-C,C]$. By {\rm (A2)}, we get that
\begin{align*}
|f(x,y,\xi)|\le\left|x\right|^{2k}\psi(x,y),\quad
|F(x,y, \xi)| \le \int\limits_{0}^{\xi}\left|x\right|^{2k}\psi(x,y) \dif \tau \le C\left|x\right|^{2k}\psi(x,y).
\end{align*}
Hence for all $(x,y)\in\Omega,\; 0\le s\le t\le C$, we have
\begin{gather*}
f(x,y, s)s-2F(x,y, s)\le f(x,y, t)t - 2F(x,y, t) +6C\left|x\right|^{2k}\psi(x,y).
\end{gather*}
Thus we get the lemma when $0\le s\le t$. Similarly we obtain it if $t\le s\le 0$.
 The proof of Lemma \ref{lemmcc0} is complete.
\end{proof}
We now show the main lemma of this paper.
\begin{lemma} \label{lemmcc1}
Assume that $f$ satisfies conditions {\rm (A1)--(A3)} and {\rm (A5)}.
Then $\Phi$ satisfies the $(C)_c $ condition for all $c\in \mathbb R$.
\end{lemma}
\begin{proof}
Let $\{u_n\}_{n=1}^{+\infty} \subset \mathscr S_{2,0 }^1(\Omega)$ be a $(C)_c$ sequence, i.e.,
\begin{align}\label{TriLuyen8: DG 23}
\Phi(u_n) \to c\; \mbox{ as }\; n \to +\infty, \quad 
\lim\limits_{n \to +\infty} \left(1+\left\| u_n\right\|_{\mathscr S_{2,0 }^1(\Omega)}\right) \left\| 
\Phi'(u_n)\right\|_{(\mathscr S_{2,0 }^1(\Omega))^*} =0, 
\end{align}
hence 
\begin{equation}
 \begin{aligned}
&\lim\limits_{n\to +\infty}\int\limits_{\Omega}\left(\frac{1}{2}f(x,y,u_{n}(x,y))u_{n}(x,y) - F(x,y,u_{n}(x,y))\right)\dif x\dif y\\
& = \lim\limits_{n\to +\infty}\left(\Phi(u_{n})-\frac{1}{2}
\langle \Phi'(u_{n}),u_{n}\rangle \right) = c.
\end{aligned}\label{lemmcc11}
\end{equation}
We first show that $\{u_n\}_{n=1}^{+\infty}$ is bounded in $\mathscr S_{2,0 }^1(\Omega)$
by a contradiction argument. By passing to a sub-sequence if necessary, we can assume that
 $\left\|u_n\right\|_{\mathscr S_{2,0 }^1(\Omega)}>1$ and
\begin{align}\label{TriLuyen8: DG 25}
\left\|u_n\right\|_{\mathscr S_{2,0 }^1(\Omega)} \to +\infty  \mbox{ as } n \to +\infty.
\end{align}
Setting
$$ 
w_{n}= \frac{u_{n}}{\|u_{n}\|_{\mathscr S_{2,0 }^1(\Omega)}},
 $$
then $\left\|w_n\right\|_{\mathscr S_{2,0 }^1(\Omega)} =1$, so we can extract a sub-sequence relabelled $\{w_{n}\}_{n=1}^{+\infty}$ such that
 $\{w_{n}\}_{n=1}^{+\infty}$ converges weakly to $w$ in $\mathscr S_{2,0 }^1(\Omega)$. Since $\Omega $ is bounded, Theorem \ref{DL:3} implies that
\begin{eqnarray} \label{TriLuyen5: DLP11}
w_n &\to& w \mbox{ strongly in } L^p_k(\Omega), 1\le p < 2_k \mbox{ as } n \to +\infty,  \notag\\
w_n &\to & w \mbox{ a.e. in } \Omega \mbox{ as } n \to +\infty. 
\end{eqnarray}
Now, we consider two possible cases: $w=0$ and $w \ne 0$
\smallskip

\noindent\textbf{Case 1:} If  $w=0$ then for any $n \in   \mathbb N$ there exists $t_n \in [0,1]$ such that
\begin{eqnarray} \label{TriLuyen5: DLP12}
\Phi(t_{n}u_{n}) = \max\left\{\Phi(su_{n}) : s \in [0,1]\right\}.
\end{eqnarray}
 Fix a positive integer $m$ and put $v_{n}= (4m)^{1/2}w_{n}$ for
every positive integer $n$. Therefore, we have
\begin{eqnarray}\label{TriLuyen5: BS11}
v_n &\rightharpoonup& 0 \mbox{ weakly in } \mathscr S_{2,0 }^1(\Omega)\ \mbox{ as } n \to +\infty, \notag \\
v_n &\to& 0  \mbox{ a.e. in  } \Omega \ \mbox{ as } n \to +\infty,  \\
v_n &\to& 0 \mbox{ strongly in } L^p_k(\Omega), 1\le p < 2_k \ \mbox{ as } n \to +\infty. \notag
\end{eqnarray}
By {\rm (A1)}, applying H\"{o}lder's  inequality, we get
\begin{align}\label{TriLuyen5: BS12}
\left| \int\limits_{\Omega}F(x,y,v_{n}(x,y))\dif x\dif y\right| \le   \int\limits_{\Omega}\left|x\right|^{2k}\left|v_n(x,y)\right|^{q_1-1}\left|v_n(x,y)\right| \dif x\dif y
+ C_0  \int\limits_{\Omega}\left|x\right|^{2k}\left|v_n(x,y)\right| \dif x\dif y \notag \\
\le \left(\int\limits_{\Omega}\left|x\right|^{2k}\left|v_n(x,y)\right|^{q_1}\dif x\dif y  \right)^{\frac{q_1-1}{q_1}}
 \left(\int\limits_{\Omega}\left|x\right|^{2k} \left|v_n(x,y)\right|^{q_1}\dif x\dif y  \right)^{\frac{1}{q_1}} + C_0  \int\limits_{\Omega}\left|x\right|^{2k}\left|v_n(x)\right| \dif x\dif y.
\end{align}
From \eqref{TriLuyen5: BS11} and \eqref{TriLuyen5: BS12}, we have
\begin{align*}
\left| \int\limits_{\Omega}F(x,y,v_{n}(x,y))\dif x\dif y\right| \to 0 \mbox{ as } n \to +\infty,
\end{align*}
hence
\begin{align*}
\lim_{n\to+\infty}\int\limits_{\Omega}F(x,y,v_{n}(x,y))\dif x\dif y =0.
\end{align*}
Since $\lim\limits_{n\to+\infty}(4m)^{1/2}\|u_{n}\|_{\mathscr S_{2,0 }^1(\Omega)}^{-1} =0$, there exists an
integer $N_{m}$ such that 
\[
\Phi(t_{n}u_{n})\ge \Phi(v_{n}) = 2m - \int\limits_{\Omega}F(x,y,v_{n}(x,y))\dif x\dif y \ge m\quad
\forall n\ge N_{m},
\]
that is, $\lim\limits_{n\to+\infty}\Phi(t_{n}u_{n})=+\infty$. Since $\Phi(0)=0$ and
$\lim\limits_{n\to+\infty}\Phi(u_{n})=c$, it implies $t_{n} \in (0,1)$ for any
sufficiently large $n$ and
\begin{align*}
\int\limits_{\Omega}|\nabla_{G_k}(t_{n}u_{n})|^{2} \dif x\dif y
- \int\limits_{\Omega}f(x, y, t_{n}u_{n}(x,y ))t_{n}u_{n}(x,y )\dif x\dif y
&= \langle \Phi'(t_{n}u_{n}),t_{n}u_{n} \rangle \\
&= t_{n}\frac{\dif}{\dif t}\Big|_{t=t_{n}}\Phi(tu_{n})=0.
\end{align*}
 Therefore, by Lemma \ref{lemmcc0}, we get
\begin{align*}
&\int\limits_{\Omega}\left(\frac{1}{2}f(x,,y, u_{n}(x,y))u_{n}(x,y)-F(x,y,u_{n}(x,y))\right)\dif x\dif y \\
&\ge \int\limits_{\Omega}\left(\frac{1}{2}f(x,y,t_{n}u_{n}(x,y))t_{n}u_{n}(x,y)-F(x,y,t_{n}u_{n}(x,y))\right)\dif x\dif y
 -C_3\|\psi\|_{L_k^1(\Omega)} \\
&= \int\limits_{\Omega}\left(\frac{1}{2}|\nabla_{G_k}(t_{n}u_{n})|^{2}-F(x,y, t_{n}u_{n}(x,y))\right)\dif x\dif y
 -C_3\|\psi\|_{L_k^1(\Omega)} \\
&= \Phi(t_{n}u_{n})-C_3\|\psi\|_{L^1_k(\Omega)}\to+\infty,
\end{align*}
which contradicts  \eqref{lemmcc11}.

\noindent\textbf{Case 2:} 
If $w\neq 0$, the Lebesgue measure of the set $\Theta  = \{(x,y)\in\Omega : w(x,y)\neq 0\}$
 is positive. We have $\lim\limits_{n\to +\infty}|u_{n}(x,y)|=+\infty$ for every $(x,y)$
in $\Theta$. Hence, by (A4) we deduce
\begin{align}\label{TriLuyen5: BS14}
 \frac{f(x,y,u_{n}(x,y))}{ u_{n}(x,y)} \to +\infty \mbox{ as  } n \to +\infty.
\end{align}
Using \eqref{TriLuyen5: BS14}, (A3) and the fact that $\Phi(u_n) \to c$,
we deduce via the  generalized Fatou lemma that
{\allowdisplaybreaks
\begin{align*}
\frac 12&= \liminf_{n\to +\infty}\left[\frac{1}{2} - \frac{\Phi(u_{n})}{\|u_{n}\|^2_{\mathscr S_{2,0 }^1(\Omega)}}\right]
 = \liminf_{n\to +\infty} \int\limits_{\Omega}\frac{F(x,y,u_{n}(x,y))}{\|u_{n}\|_{\mathscr S_{2,0 }^1(\Omega)}^{2}}\dif x\dif y \\
&= \liminf_{n\to+\infty}\int\limits_{\Omega}\int\limits_{0}^{1}
 \frac{f(x,y,\xi u_{n}(x,y))}{\xi u_{n}(x,y)} \cdot \xi |w_{n}(x,y)|^{2} \,\dif \xi\, \dif x\dif y \quad  (\mbox{since  } w\equiv 0 \mbox{ on  } \Omega\backslash \Theta)\\
&= \liminf_{n\to+\infty}\int\limits_{\Theta}\int\limits_{0}^{1}
 \frac{f(x,y,\xi u_{n}(x,y))}{\xi u_{n}(x,y)}\cdot \xi |w_{n}(x,y)|^{2} \,\dif \xi\, \dif x\dif y\\
&\ge \int\limits_{\Theta}\int\limits_{0}^{1}\liminf_{n\to+\infty}
 \frac{f(x,y,\xi u_{n}(x,y))}{\xi u_{n}(x,y)}\cdot \xi |w_{n}(x,y)|^{2}\,\dif \xi\,\dif x\dif y
=+\infty,
\end{align*}}
which is impossible. In any case, we obtain a contradiction.
Thus, the sequence $\{u_{n}\}_{n=1}^{+\infty}$ is bounded in $\mathscr S_{2,0 }^1(\Omega)$.
Therefore we can (by passing to a sub-sequence 
if necessary) suppose that 
\begin{eqnarray}\label{TriLuyen5: BS1}
u_n &\rightharpoonup& u \mbox{ weakly in } \mathscr S_{2,0 }^1(\Omega)\ \mbox{ as } n \to +\infty, \notag \\
u_n &\to& u  \mbox{ a.e. in  } \Omega \ \mbox{ as } n \to +\infty,  \\
u_n &\to& u \mbox{ strongly in } L^p_k(\Omega), 1\le p < 2_k\ \mbox{ as } n \to +\infty. \notag
\end{eqnarray}
Thus by {\rm (A1)}, we have
\begin{align}
&\left|\int\limits_\Omega f(x,y,u_n(x,y))(u_n(x,y) -u(x,y)) \dif x\dif y \right| \notag \\
  \le &C_0 \int\limits_\Omega \left|x\right|^{2k} \left|u_n(x,y) -u(x,y)\right|\dif x\dif y+
 \int\limits_\Omega \left|x\right|^{2k} \left|u_n(x,y) -u(x,y)\right| \left|u_n(x,y)\right|^{q_1-1}\dif x\dif y \\
\le  & C_0 \left\|u_n -u\right\|_{L_k^1(\Omega)}+
 \left(\int\limits_\Omega\left|x\right|^{2k} \left|u_n(x,y) -u(x,y)\right|^{q_1} \dif x\dif y\right)^{\frac{1}{q_1}} 
\left( \int\limits_\Omega  \left|x\right|^{2k}\left|u_n(x,y)\right|^{q_1}\dif x\dif y\right)^{\frac{q_1-1}{q_1}}.  \notag 
\end{align}
In view of  \eqref{TriLuyen5: BS1}, we can conclude that 
$$ \int\limits_\Omega f(x,y,u_n(x,y))(u_n(x,y) -u(x,y)) \dif x\dif y \to 0 \mbox{ as } n \to +\infty. $$
Thus 
\begin{align}\label{TriLuyen5: BS2}
\int\limits_\Omega  \left(f(x,y,u_n(x,y)) - f(x,y, u(x,y))\right) (u_n(x,y) -u(x,y)) \dif x\dif y \to 0 \mbox{ as } n \to +\infty. 
\end{align}
It follows from $\lim\limits_{n \to +\infty} \Phi'(u_n)=0$ and \eqref{TriLuyen5: BS1} that
\begin{align}\label{TriLuyen5: BS4}
\left \langle \Phi'(u_n) - \Phi'(u) , u_n -u \right\rangle \to 0 \mbox{ as } n \to +\infty. 
\end{align}
By \eqref{TriLuyen5: BS2} and \eqref{TriLuyen5: BS4}, we obtain
$$  \int\limits_\Omega\left| \nabla_{G_k} u_n -\nabla_{G_k} u\right|^2\dif x\dif y  \to 0 \mbox{ as } n \to +\infty.  $$
Therefore, we conclude that $u_n \to u$ strongly in $ \mathscr S_{2,0 }^1(\Omega)$.  The proof of Lemma \ref{lemmcc1}
 is complete.
\end{proof}

\begin{proof}[Proof of  Theorem \ref{thm1}]
By Lemmas \ref{lemma2}, \ref{lemma1} and  \ref{lemmcc1} all conditions of  Lemma \ref{Luyen:BDP} are satisfied. 
Thus, the problem \eqref{Luyen2}-\eqref{Luyen3} has a nontrivial weak solution.
\end{proof}

\end{document}